# An Improved Algorithm based on Shannon-Happ Formula for Calculating Transfer Function from Signal Flow Graph and Its Visualization

Hongyu Lu [1]    Chongguang Wu [2]     Shanglian Bao [1]

([1] Beijing City Key Laboratory of Medical Physics and Engineering, Peking University, China; [2] College of Information Science & Technology, Beijing University of Chemical Technology, Beijing 100029, China)

**Abstract:** A new method based on Shannon-Happ formula to calculate transfer function from Signal Flow Graph (SFG) is presented. The algorithm provides an explicit approach to get the transfer function in a format with both numerical and symbolic expressions. The adoption of the symbolic variable in SFG, which could represent the nonlinear item or the independent sub-system, is achieved by variable separation approach. An investigation is given for the solutions of several special conditions of SFG. To improve the efficiency of the algorithm, a new technique combined with Johnson method for generating the combinations of the non-touching loops is developed. It uses the previous combinations in lower order to get the ones in higher order. There is an introduction about the visualization of SFG and the subroutines for system performance analysis in the software, AVANT.

## 1  Introduction

A complicated linear system can be expressed in the format of Signal Flow Graph (SFG) concisely and intuitionisticly. The SFG takes the advantage of getting its form directly from the block diagrams. Based on SFG, it can obtain the algebraic expression of the transfer function, which is very important to the system analysis and design. Currently, the methods to calculate the transfer function with SFG are, direct elimination method, Mason's rule and the approach of solving linear algebraic equations, etc.

In electrical engineering, using the Mason's rule is a normal approach, but it demands the removal of the loops touching the $k$-th forward path [1]. The Shannon-Happ formula [2] gives an explicit approach by creating a closed SFG. It only requires finding all the loops in SFG, and has no need to concern about the forward path. Hence, The algorithm presented is developed by solving the Shannon-Happ formula.

Based on the Johnson method, a scheme is designed for obtaining the total combinations of non-touching loops. The high order loop combinations are generated step by step through the pervious combination in lower order. This approach accelerates the computing speed. Further, analysis is given for three cases required for special treating in SDG.

MacNamee introduced the expression of the transfer function all in format of symbolic variables by using the unilateral graph techniques and Kirchhoof's law [3]. Through introducing the technique of variable separation, a method capable of calculation of the transfer function from SFG with both numerical and symbolic approaches is developed. This approach takes less efforts in algorithm design. The symbols can represent the constants, nonlinear items and sub-systems so that it expand the applicability of SFG. Finally, the visualized SFG software for system simulation and performance analysis, AVNT, is presented with its structure and functions.

## 2    The Shannon-Happ formula

The reduction of the Mason's equation can directly result an expression of the transfer function. This rule must require the SFG has a source or input node, which has only branches going out, and a sink or output node, which has only branches going in as shown in Fig 1. This kind of graph is named as opened SFG. The Mason formula can be defined as,

$$G = \frac{1}{\Delta} \sum_i P_i \Delta_i \tag{1}$$

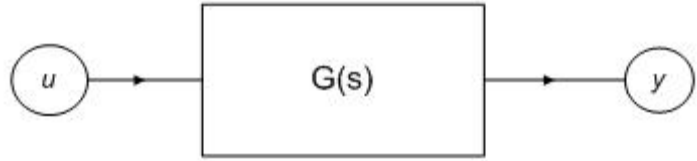

Figure 1.    An opened SFG

The meaning of the each part in Eq. 1 can refer to [1]. The difficulty of computing the $\Delta_i$ related to the forward path in Mason's rule makes structure of the algorithm complicated.

Shannon and Happ invented the closed SFG, which means to from a special loop including the input , the out nodes, and the block of $G(s)$. This task can be achieved by adding a branch with the gain of $1/G(s)$ from the output node to the input node as illustrated in Fig. 2.

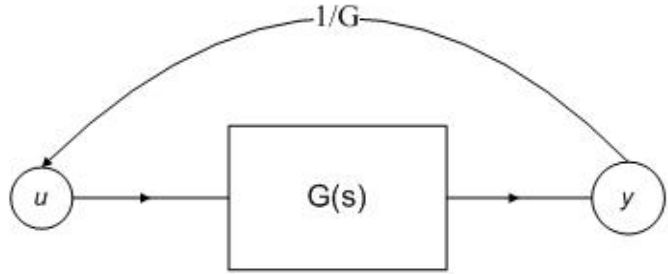

Figure 2.   A closed signal flow graph by adding a $1/G$ branch from the output node to the input node

Shannon-Happ formula can be expressed as,

$$f(1/G) = 1 - \sum_j g_j + \sum_{m,n} g_m g_n - \sum_{p,q,r} g_p g_q g_r + \cdots = 0 \tag{2}$$

Where $g_i$ has the same meaning as that in the $\Delta$ of Mason's rule.

It should be noticed that the Eq. 2 is actually a function of $1/G(s)$. The transfer function $G(s) = B(s)/A(s)$ can be gotten through,

$$G(s) \times 1/G(s) = \frac{B(s)}{A(s)} \times 1/G(s) = 1 \tag{3}$$

Thus, it has,

$$A(s) - B(s) \times 1/G(s) = 0 \tag{4}$$

The Eq. 4 is a linear equation of $1/G$, an alterative format of Eq. 2. With Shannon-Happ formula, it only needs to compute the product of the gains of the all combinations of non-touching loops. Therefore, it can greatly simplify the implementation of the algorithm.

Based on this design, the task for searching all the non-touching loops can be achieved by Johnson method [4]. Johnson method searches the fundamental loops in a directed graph, which has high efficiency and simple data structure.

Johnson method can find the self-loop automatically. There is no problem for getting the gain of the loops like $\{n,n\}$, where $n$ is the node number. This method can not process the condition that one node has two or more branches with same direction connected to another node. This kind of graph is illustrated in Fig. 3 (a). The problem can be solved by adding an additional node to link the two nodes, for example, the node No. 3 as that in Fig 3 (b).

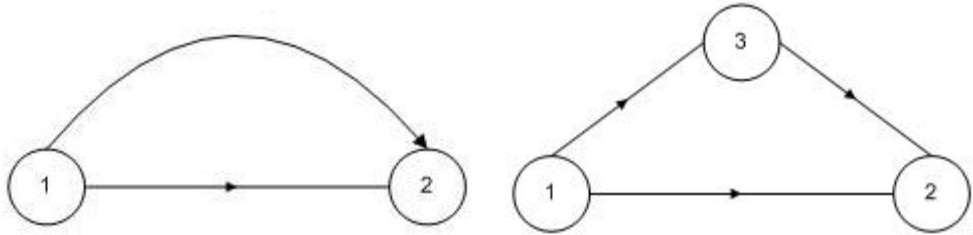

(a) Original signal flow graph      (b) Signal flow graph after inserting a new node

Figure 3. Solution to condition of two nodes with two directed branches

Another question for computing transfer function with SFG is that the input node has branches going in, or the output node has branches going out. In this case, adding the new input or output node with branch gain $1.0$ is the only solution.

## 3  Generation of the Combinations of Non-touching Loops

Computing the combination of non-touching loops is normally a task of heavy duty if always checking whether the two loops are contacted pair by pair. The number of four pair-wise non-touching loops is less than number of three non-touching loops. Moreover, in a set of combination of four pair-wise non-touching loops, any three of them should not touch each other. If treating this question in a reversed view, it is possible to add a special loop to the set of three non-touching loops to form the set of four non-touching loops. It is convenient to define the $n$ pair-wise non-touching loops as the non-touching loops in the $n$-th order, or in abbreviation, the $n$-th order loops. Supposing if there is a combination of non-touching loops in the $i$-th ($i > 2$) order, the total number of combinations are $m$, which is shown in table 1.

Table 1. The $m$ combination of non-touching loops in the $i$-th order, $i > 2$

$$\begin{matrix} l_{11} & l_{12} & \cdots & \cdots & l_{1i} \\ l_{21} & l_{22} & \cdots & \cdots & l_{2i} \\ \vdots & \vdots & \ddots & & \vdots \\ \vdots & \vdots & & \ddots & \vdots \\ l_{m1} & l_{m2} & \cdots & \cdots & l_{mi} \end{matrix}$$

In Table 1, the $l_{jk}$ corresponds to the loop number. It has the feature, $l_{j1} \geq l_{j2} \cdots \geq l_{ji}$., which $l_{j,k} = l_{j,k+1}$ means a self-loop. Further, it demands elements in the first column having the attribute, $l_{j+1,1} \geq l_{j,1}$. There are $p$ different loops in the first column. It is convenient to build the non-touching loops in the 2nd order first with an independent subroutine.

If adding a new loop to the $i$-th order loops to form the $(i+1)$-th order loops, the new loops combination should abide those two ranking features, too. The $p$ different loops, $l_{1,1} \geq l_{2,1} \cdots \geq l_{m,1}$ in the first column is arranged as $L_1, L_2, \cdots, L_{p-1}, L_p$, which has the feature, $L_k < L_{k+1}$. It can approve that new loops able to put into the $j$ row in Table 1 are all in the set, $A = \{L_1, L_2, \cdots, L_{p-1}\}$. Taking a $L_k$ from $A$, which has $L_k > l_{j1}$, and insert it to the first column of every possible row in Table 1, can create a new non-touching loop in the $(i+1)$-th order. Here, it demands detecting the $L_k$ does not contact with any exiting loops in a row. Repeating this procedure to other elements in $A$, all combinations of non-touching loops in the $(i+1)$-th can be easily gotten. By the similar way, the product of the loop gains in the $(i+1)$-th order can be gotten from the results in the $i$-th order.

## 4   Symbolic Variables

In practice, it needs to know the influence of a parameter to the system performance, especially its role in the transfer function. Many industrial process systems have nonlinear units, for example, the term for time varying delay, which its Laplace transformation is $e^{-\tau s}$. It is inevitably induce the error if utilizing the linear approximation such as Pade approximation. It is meaningful to study the effects of sub-system in the whole system, and encapsulate the sum-system as a model in the CAD software.

Treating the parametric, the nonlinear terms, and the comparably independent sub-system as symbolic variables can give a satisfactory solution to express them in transfer function. Further, applying operation such as the differential to the equations with symbolic term can get more features of the system.

Assuming there is a symbolic variable $V(s)$ in SFG, the key point is that the loops containing the branches of $V(s)$ must contact each other. Hence, only the 1st order $V(s)$ can appear in the expression of $G(s)$. It is not possible for $G(s)$ to have high order like $V^k(s)$ $(k \geq 2)$ term.

The numerical Shannon-Happ formula is not a normal polynomial of $s$ but also contains a variable of $1/G$. This variable can be extracted through a special data matrix, which has an index table to pointing to a loop containing $1/G$ in the table1. Therefore, it can easily find the $1/G$ terms in the product of the non-touching loops gains.

The separation of one variable $V(s)$ can take the same step by creating a new data matrix with index table in the program. More symbolic variables can be applied in SFG by the same way. The implementation method relates to the SFG structure, the data matrices, and the explanation of the separated variables. It should be point out that this approach will make the expression of $G(s)$ very complicate for analysis if there are too many variables in SFG.

## 5  The Procedure for Calculating Transfer Function with Symbolic Variables

The implementation of solving the function of Shannon-Happ formula with one symbolic variable example can be described as,

1) Read the topological network information of SFG, and the branch-parameter matrix including branches of $1/G$ and $V(s)$ Then Johnson method is applied to get the basic loops.

2) Adopt the standard combination-making algorithm to generate the 2nd order non-touching loops. Use the method in Section 4 for generating the non-touching loops to create all the $m$ combinations of the $i$-th ($i \geq 3$) order non-touching loops.

3) Compute the gain of each loop. Then search in the index table to find whether the gain has the term $1/G$ or $V(s)$. Hence, a loop gain can be expressed explicitly by the constant $c$ and the its $s^k$ in Laplace domain, and the term of $1/G$ or $V(s)$. They are in four kinds of format, $cs^k$, $cs^k \cdot V(s)$, $cs^k \cdot 1/G$, and $cs^k \cdot V(s) \cdot (1/G)$. The gain can define as $g_{j\_i}$.

4) Calculate the product of each row of Table 1. For example, to the row $j$, there has $g_j = g_{j\_1} g_{j\_2} \cdots g_{j\_i}$. Based on the variable separation method, the $g_i$ still takes the four kinds of format as in step 3.

5) Make the sum, $G_j = \sum_{j=1}^{m} g_j = \sum_{j=1}^{m} g_{j\_1} g_{j\_2} \cdots g_{j\_i}$, according to condition of the power of $s$, and whether it has the $1/G$ or $V(s)$. Then $G_i$ is normalized in order.

6) Calculate Shannon-Happ formula with $f(1/G) = 1 - G_1 + G_2 - G_3 + \cdots$. This step extracts the variables, and normalizes the form of $f(1/G)$ according to the power of of $s$, $1/G$ or $V(s)$. In the result, the polynomial with $1/G$ is the numerator, and the polynomial without is $1/G$ is the denominator. Hence, the final normalized algebraic expression of $G(s)$ is found.

## 6 Discussion

This algorithm can be applied to the discrete system. It only needs to replace $s$ with $z$. To a MISO system with zero status, according to the principle of linear superposition, it has,

$$Y_f(s) = G_1(s)F_1(s) + G_2(s)F_2(s) + \cdots\cdots + G_n(s)F_n(s) \tag{5}$$

To a SISO system with the $n$ non-zero status system, it has,

$$Y(s) = Y_X(s) + Y_f(s) = \frac{1}{A(s)}\sum_i^n M_i(s)X_i(0) + \frac{B(s)}{A(s)}F(s) \tag{6}$$

Where $Y_X(s)$ denotes the zero input responses, and $X_i(0)$ denotes the initial status of the component $i$ in Laplace domain. The $X_i(0)$ is treated as the input of the system. Then the output is $Y_{X\_i}(s)$. Hence, the sum of all the components with initial status is $Y_X(s) = \sum_{i=1}^m Y_{X\_i}(s)$. To the MIMO system, it has the similar expression under the principle of linear superposition.

## 7 Introduction to AVANT Software

AVANT software is a visual simulation platform with plenty of features based on the proposed algorithm for computing the transfer function in SFG. The software structure of AVANT is illustrated in Fig. 4.

The visual graphics editor defines the graphic icons, and their manipulations complying with the widely used visual programming software. It provides the node, directed polygonal line, self-loop, and directed arc with their corresponding value labels. The dynamic operations to those graphics icons include creation, display, management, selection, movement and removal of the components. It can identify the attributes and the relationships of the linkage of these graphics icons. The manipulations of a graphic icon can response with the actions of the related icons. The topological and branch parameters of the SFG can be built with information of the graphic icons on the visual graphics editor. Drawing the directed line or arc with $1/G$ gain between two nodes assigns the output node and input node for SFG automatically. In addition, the visual editor calls the subroutine of Shannon-Happ formula for computing the transfer function after finishing the graphics design.

The normalized algebraic expression of the transfer function supplies a powerful tool for system performance analysis. The models in AVANT include, using the continued fraction approach to give a transfer function in the reduced order, drawing the system frequency response curve and the Nyquist diagram, analyzing the stability by Ruth's rule, adoption of QR method to find the system poles and zeros.

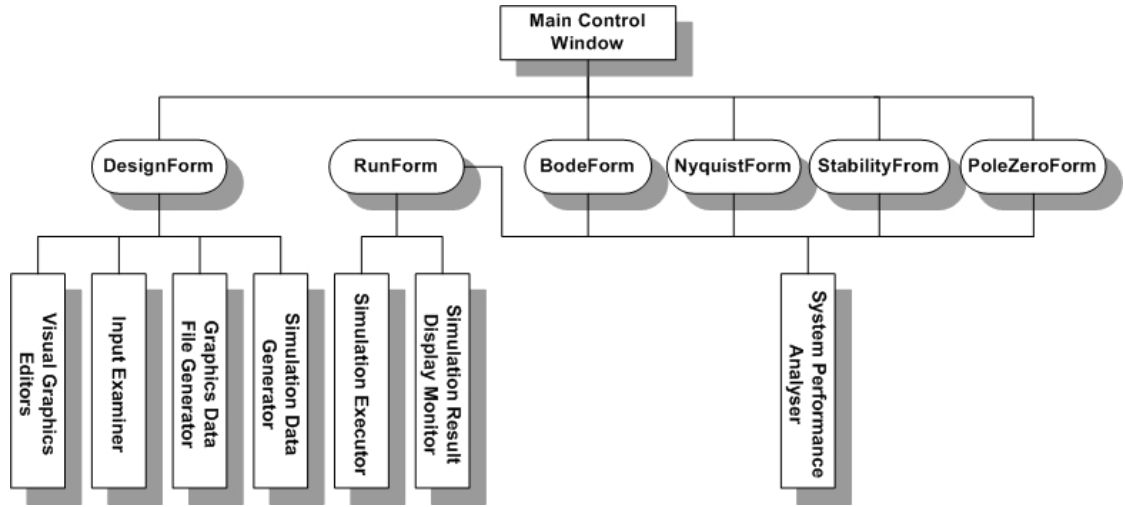

Figure 4. The software structure of AVANT

## 8  Experiments and Results

Fig. 5 shows a SFG with 7 nodes, in which the node 1 is inserted to get an input node on the original SFG. The node 3 to node 4 is a special case so that the algorithm inserts a new node, No. 8, on the branches with the gain of 0.5. If assuming $G(s) = B(s)/A(s)$, in which $B(s)$ is the numerator, and $A(s)$ is the denominator, the proposed algorithm gives the results of transfer function as illustrated in Table 2.

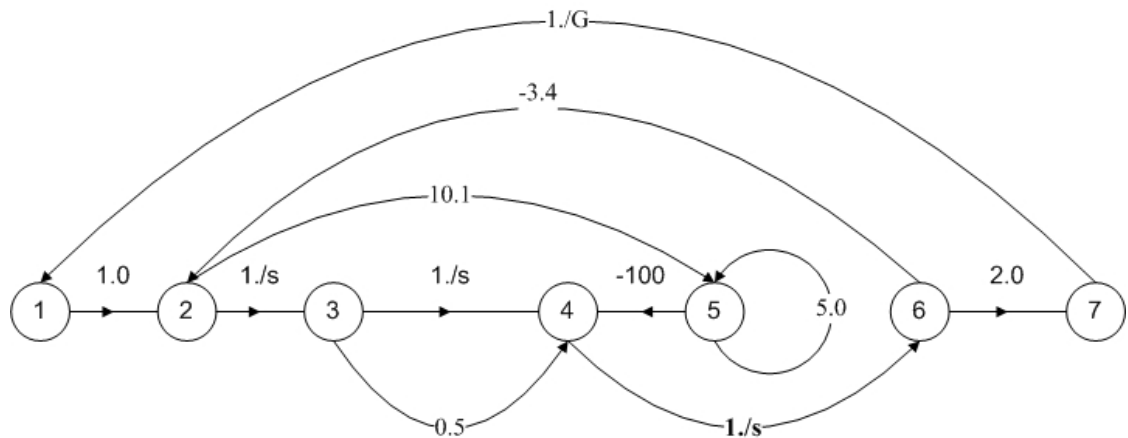

Figure 5.  The signal flow graph with 7 nodes including a self-loop

Table 2. The comparison for real value and result of the program for the transfer function

| Power of s | Coefficients in B(s) | | Coefficients in A(s) | |
|---|---|---|---|---|
| | Real value | Result of program | Real value | Result of program |
| 0 | 3.4 | 3.400000 | 2 | 2.000000 |
| 1 | 1.7 | 1.700000 | 1 | 1.000000 |
| 2 | 858.5 | 858.500122 | 505 | 505.000031 |
| 3 | 1.0 | 1.000000 | 0 | 0.000000 |

Fig. 6 shows a SFG with a sub-system includes a time delay component, which can be represented by a symbolic variable $V(s)$. The transfer function can be calculated as,

$$G(s) = \frac{1}{s+1} \times \frac{s+4}{s+2} \times \frac{e^{-0.67s}}{s+3} \times 2 = \frac{2s+8}{s^3+6s^2+11s+6} e^{-0.67s} \quad (7)$$

If treating the sub-system, which is surrounded by a dashed block, as a variable $V(s)$, it can get $V(s) = e^{-0.67s}/(s+3)$. The result is shown in Table 3.

Table 3. The result of calculation the transfer function from SFG with symbolic variable $V(s)$

| Power of s | Coefficients in B(s) | | Coefficients in A(s) | |
|---|---|---|---|---|
| | Without V(s) | With V(s) | Without V(s) | With V(s) |
| -2 | 0.0 | 8.0 | 2.0 | 0.0 |
| -1 | 0.0 | 2.0 | 3.0 | 0.0 |
| 0 | 0.0 | 0.0 | 1.0 | 0.0 |

From Table 3, it can get

$$G(s) = \frac{2sV(s) + 8V(s)}{s^2+3s+2} = \frac{2s+8}{s^2+3s+2} \times \frac{e^{-0.67s}}{s+3} = \frac{2s+8}{s^3+6s^2+11s+6} e^{-0.67s} \quad (9)$$

Hence, the result of $G(s)$ in Eq. 9 is as same as Eq. 10.

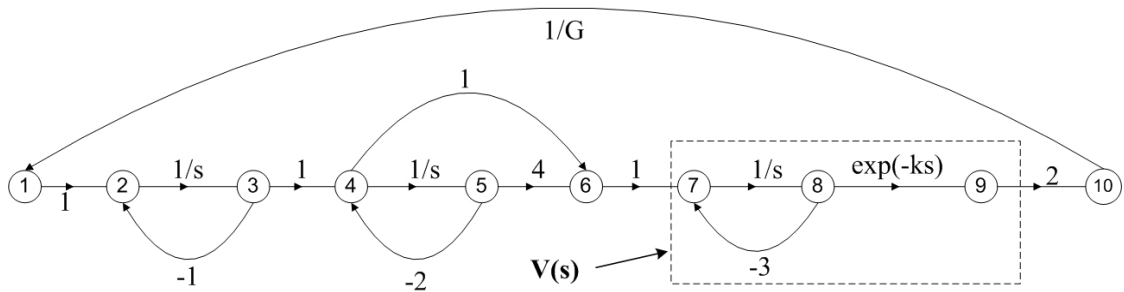

Figure 6.  The signal flow graph with 10 nodes including a sub-system with symbolic variable $V(s)$

Fig. 7 gives an example of a visualized SFG of a CSTR control system by using the visual graphics editor of AVANT. It needs to get the transfer function from node 1 to node 13. The reason to insert node 18 is to from an output node. Introducing the $1/G$ branch from node 18 to node 1 is to build a closed SFG for utilizing the Shannon-Happ formula. Fig. 8 shows the simulation result of $G(s)$ and its 3rd reduced order expression. Table 4 gives the comparison of the result by calling the proposed algorithm, and the result in the reference [5]. Fig. 9 displays the two frequency response curves from the original transfer function and its 3rd order expression. It can be seen that the two response curves are nearly superposed. Therefore, from the view of frequency response, the channel from node 1 to node 13 can be simplified to a simple 3rd order system.

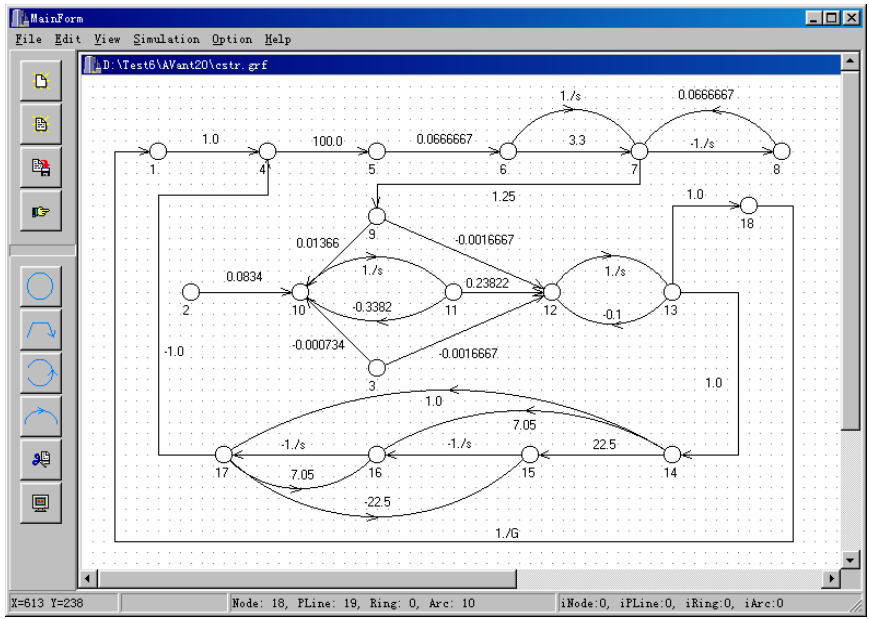

Figure 7. The Visual signal flow graph of a CSTR control system

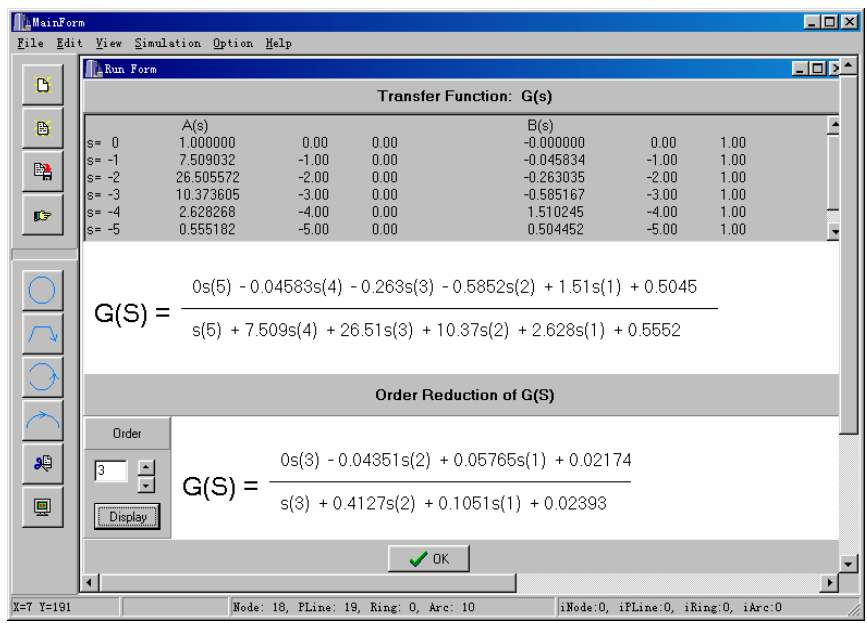

Figure 8. The two transfer functions of the CSTR control system

Table 4. The comparison of the transfer function from node 1 to node 13

| Power of s | Coefficients in B(s) | | Coefficients in A(s) | |
| --- | --- | --- | --- | --- |
| | Value in reference [10] | Result of program | Value in reference [10] | Result of program |
| 0 | 0.50445086 | 0.504452 | 0.55518079 | 0.555182 |
| 1 | 1.5102396 | 1.510245 | 2.6282597 | 2.628268 |
| 2 | -0.58516490 | -0.585167 | 10.373603 | 10.373605 |
| 3 | -0.26303399 | -0.263035 | 26.505478 | 26.505572 |
| 4 | -0.045834191 | -0.045834 | 7.5090303 | 7.509032 |
| 5 | 0.0 | 0.000000 | 1.0000000 | 1.000000 |

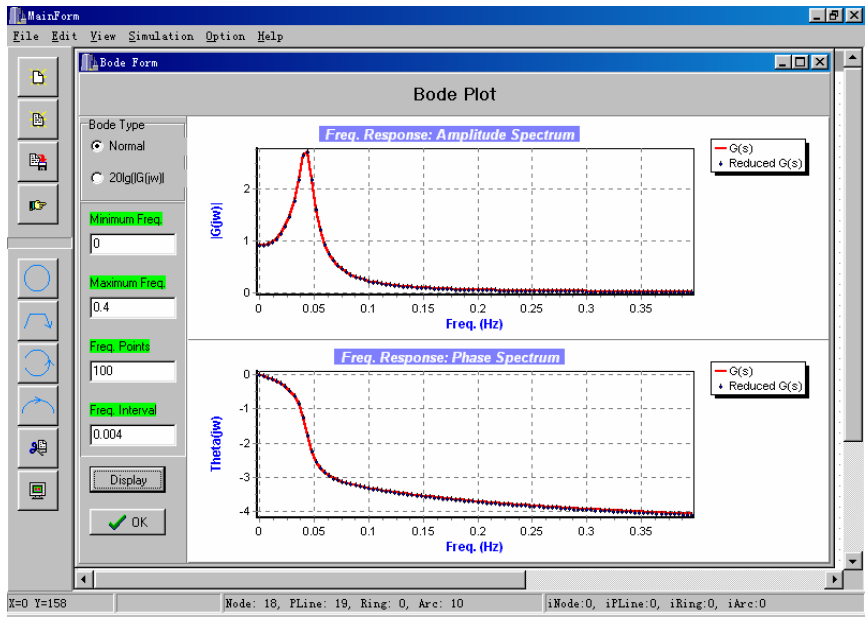

Figure 9. The two frequency repose curves of $G(s)$, the original one and its 3rd reduced order

## 9 Conclusion

From the introduction about the algorithm and demonstrations of experiments, the improved method based on Shannon-Happ formula contributes several concise and fast techniques for finding the normalized algebraic expression of the transfer function in SFG. The inducing of symbolic variables makes this technology adapt to more general and practical fields. The visual platform, AVANT, provides the strong abilities to create SFG, to calculate the transfer function, and perform system analysis.